\documentclass{amsart}
\usepackage{amsmath,amssymb}
\newtheorem{theorem}{Theorem}
\newtheorem{lemma}[theorem]{Lemma}
\newtheorem{proposition}[theorem]{Proposition}
\newtheorem{corollary}[theorem]{Corollary}

\begin{document}

\title[Local rings of
  bounded module type]{Local rings of
  bounded module type are almost maximal valuation rings}
\author{Fran\c{c}ois Couchot}

\begin{abstract} It is shown that every commutative
local ring of bounded module type is an almost maximal valuation ring.
\end{abstract}
\maketitle

We say that an associative commutative unitary ring $R$ is of bounded
module type if there exists a positive integer $n$ such that every
finitely generated $R$-module is a direct sum of submodules generated by
at most $n$ elements. For instance, every Dedekind domain has bounded
module type with bound $n=2$. Warfield \cite[Theorem 2]{War} proved
that every commutative local ring of bounded module type is a valuation
ring. Moreover, Gill \cite{Gil} and Lafon \cite{Laf} showed
independently that a valuation ring $R$ is almost maximal if and only
if every finitely generated $R$-module is a direct sum of cyclic
modules. So every almost maximal valuation ring has bounded module
type and V\'amos proposed the following conjecture, in the Udine
Conference on Abelian Groups and Modules, held in 1984:

 ``A local ring of bounded module type is an almost
maximal valuation ring."

P. Zanardo \cite{Zan} and P. V\' amos \cite{Vam}  first investigated
this conjecture and proved it respectively for strongly discrete
valuation domains and $\mathbb{Q}$-algebra valuation domains. More
recently the author proved the following theorem (see \cite[Corollary
9 and Theorem 10]{Cou}).

\begin{theorem} \label{T:main} Let $R$ be a local ring of bounded
  module type. Suppose that $R$ satisfies one of the following conditions:
\begin{enumerate}
\item There exists a nonmaximal prime ideal $J$ such that 
$R/J$ is almost maximal.
\item The maximal ideal of $R$ is the union of all nonmaximal prime ideals.
\end{enumerate}
Then $R$ is an almost maximal valuation ring.
\end{theorem}

From this theorem we deduce that it is enough to show the following 
proposition to prove V\'amos' conjecture : see \cite[remark
1.1]{Cou}. Recall that a valuation ring is \textit{archimedean} if the maximal
ideal is the only non-zero prime ideal.

\begin{proposition} \label{P:main}  Let $R$ be an  archimedean valuation ring
  for which there exists a positive integer $n$ such that every
finitely generated uniform module is generated by at most $n$
elements. Then $R$ is almost maximal. 
\end{proposition}

 In the sequel we use the same terminology and the same
notations as in \cite{Cou} except for the symbol
$A\subset B$ which means that $A$ is a proper subset of $B$. We also
use the terminology of \cite{Cou1}. Some results of \cite{Cou1} and the
following lemmas will be useful to show proposition~\ref{P:main}.  

An $R$-module $E$ is said to be \textit{fp-injective}(or
absolutely pure) if $\mathrm{Ext}_R^1(F,E)=0,$ for every finitely
  presented $R$-module $F.$ A ring $R$ is called \textit{self
    fp-injective} if it is fp-injective as $R$-module.

\begin{lemma} \label{L:inj} Let $R$  be a coherent archimedean
  valuation ring,  $Y$ an injective
$R$-module, $X$ a pure submodule of $Y$ and $Z = Y/X$. Then $Z$ is injective.
\end{lemma}
\emph{Proof.}~ We have $X$ fp-injective because it is a pure
submodule of an injective module. So,  since $R$ is coherent, $Z$ is
also $fp$-injective by
\cite[Theorem 1.5]{Cou2}. Let
$J$ be an ideal of $R$ and $f : J\longrightarrow Z$ an homomorphism. By
\cite[corollary 36]{Cou1}, $J$ is countably
generated. Then there exists a sequence ${(a_n)}_{n\in\mathbb{N}}$ of 
elements of $J$ such that $J
= \cup_{n\in\mathbb{N}} Ra_n$ and $a_n\in Ra_{n +1}\ \forall
n\in\mathbb{N}$. Since $Z$ is $fp$-injective, for every
$n\in\mathbb{N}$ there exist $z_n\in Z$
such that $f(a_n) = a_n z_n$. It follows that $a_n(z_{n+1} - z_n)=0$,
$\forall n\in\mathbb{N}$. By induction on $n$ we build a sequence
${(y_n)}_{n\in\mathbb{N}}$ of  elements of $Y$
such that  $z_n = y_n + X$ and
$a_n(y_{n+1} - y_n)=0, \ \forall n\in\mathbb{N}$. Suppose 
$y_0,\ldots,y_n$ are built. Let
$y'_{n+1} \in Y$ such that $z_{n+1} = y'_{n+1} + X$. Then $a_n 
(y'_{n+1} - y_n) \in X$.
Since $X$ is a pure submodule of $Y$ there exists $x_{n+1}\in X$ such 
that $a_n(y'_{n+1} -
y_n) = a_n x_{n+1}$. We put $y_{n+1} = y'_{n+1} - x_{n+1}$.  The 
injectivity of $Y$ implies
that there exists $y\in Y$ such that $a_n(y-y_n)=0$, 
$\forall n\in\mathbb{N}$. We set $z=
y+X$. Then it is easy to check that $f(a) = az$ for every $a\in J$. \qed

\bigskip
We say that a ring $R$ is an \textit{IF-ring} if it is coherent and
self fp-injective. Equivalent conditions for a valuation ring to be an
IF-ring are given in \cite[Theorem 10]{Cou1}.

In the four following lemmas $R$ is a non-noetherian
archimedean valuation ring and an IF-ring. In the sequel $P$ is the
maximal ideal of $R$. 
\begin{lemma} \label{L:arIF}  Then:
\begin{enumerate}
\item $P$ is the only prime ideal of $R$.
\item $P$ is not finitely generated.
\item $P$ is faithful.
\end{enumerate}
\end{lemma}
\emph{Proof.} Since $R$ is self
fp-injective, $\forall a\in P,\ (0:a)\ne 0$ by \cite[Proposition
1(1)]{Cou}. Hence $P$ is the only prime ideal.

If $P$ is finitely generated then all prime ideals of $R$
are finitely generated. It follows that $R$ is noetherian.

If $P$ is
not faithful then there exists $a\in R$ such that $P=(0:a)$. We deduce
that $P$ is finitely generated because $R$ is coherent, whence a
contradiction. So $P$ is faithful. \qed

\begin{lemma} \label{L:unifaith} Let $U$ be a non-zero uniserial
  $R$-module such that, for each $x\in U\setminus \{0\}$ $(0 : x)$ is
  not finitely generated. Then $U$ is $fp$-injective if and only if
  $U$ is faithful. 
\end{lemma}
\emph{Proof.}~ Assume that $U$ is faithful. Let $s\in R$ and 
$x\in U$  such that
$(0 : s) \subseteq (0 : x)$. We will show that $x\in sU$. Since $R$ is
coherent $(0 : s)$ is principal, whence
we have $(0 : s)\subset (0 : x)$. Hence $s(0 : x) \not= 0$. So,
$\exists u\in U$ such that $(0 : u)\subset s(0 : x)
\subseteq (0 : x)$ because $U$ is faithful and uniserial. Hence $x\in
Ru$ and there exists $r\in R$ such that $x = ru$. By \cite[Lemma
2]{Cou1} $(0 : u) = r(0 : x)$. It follows that $r(0 : x)\subset s(0 : x)$. Therefore $r\in Rs$. We get that $x = sau$ for some $a\in R$.

Conversely, assume that $U$ is $fp$-injective. Since $P$ is faithful
and the only prime ideal by lemma~\ref{L:arIF}, we can apply
\cite[Proposition 27]{Cou1}. Hence $U$ is faithful. \qed

\begin{lemma} \label{L:puruni} Let $E$ be a non-zero $fp$-injective
  $R$-module. Then $E$ contains a faithful uniserial pure submodule
  $U$. Moreover, if there exists $y\in E$ such that $(0:y)=0$ then we
  can choose $U=Ry$. Else, $U$ is not finitely generated.
\end{lemma}
\emph{Proof.}   If there exists $y\in E$ such that $(0:y)=0$, we
set $U = Ry$. Since $U\simeq R$, $U$ is fp-injective. So it is a pure
submodule of $E$. 

Suppose now that $(0 : x) \not=0, \forall x\in E$. We put $I=(0:x)$
for a non-zero element $x$ of $E$ and $J=(0:I)$. If $0\ne y\in
E(R/I)$ and $(0:y)\subseteq I$ then $(0:y)=sI$ for some $s\in
R\setminus J$ (see
\cite[Proposition IX.2.1]{FuSa}). Since $P$ is faithful and the only
prime ideal by lemma~\ref{L:arIF}, then
$E(R/I)$ is faithful by \cite[Proposition 27]{Cou1}. So $\cap_{s\in
  R\setminus J}sI=0$. On the other hand, $R$ is countably
cogenerated by \cite[Corollary 35]{Cou1}. By \cite[Lemma 30]{Cou1}
there exists a countable family $(s_n)_{n\in\mathbb{N}}$ of elements
of $R\setminus J$ such that
$\cap_{n\in\mathbb{N}}s_nI=0$ and $s_n\notin Rs_{n+1}$ for each
  $n\in\mathbb{N}$. We put $t_0=1$, $t_1=s_1$. For each integer $n\geq
  1$, let $t_{n+1}\in P$ such that $s_{n+1}=t_{n+1}s_n$. We get
  a family $(t_n)_{n\in\mathbb{N}}$ of elements 
of $R$ such that $t_0\ldots t_n I \not= 0,\ \forall n\in\mathbb{N}$,
and $\cap_{n\geq 0} t_0\ldots t_n I = 0$. By induction on $n$, we 
build a sequence $(x_n)_{n\in\mathbb{N}}$ of elements of $E$ such
that $(0:x_n)=t_0\ldots t_nI,\ \forall n\in\mathbb{N}$. We
set  $x_0 = x$. Since $t_{n+1}(0:x_n)\ne 0$ we have
$(0:t_{n+1})\subset (0:x_n)$. So, the fp-injectivity of $E$ implies
that $\exists x_{n+1}\in E$ such that $x_n =t_{n+1} x_{n+1}$. By
\cite[Lemma 2]{Cou1} it follows that
$(0:x_{n+1})=t_{n+1}(0:x_n)=t_0\ldots t_{n+1}I$. Let $U$ be the
submodule of $E$ generated by $\{ x_n \mid n\geq 0\}$. Then $U$ is
uniserial and faithful. Suppose $\exists u\in U$ such that
$(0:u)$ is principal. Since $R$ is an IF-ring there exists $s\in R$
such that $(0:u)=(0:s)$. The fp-injectivity of $E$ implies that $u=se$
for some $e\in E$. It follows that $(0:e)=s(0:u)=0$ by \cite[Lemma
2]{Cou1}. We get a contradiction. Thus we can apply
lemma~\ref{L:unifaith}. Hence $U$ is
fp-injective. So, it is a pure submodule of $E$. Then $U$ is not
finitely generated, else $U=Ry$ with $(0:y)=0$.  \qed 

\bigskip
If $N$ is a finitely generated $R$-module, $\mu(N)$ denotes its
minimal number of generators. Then $\mu(N)=\mathrm{dim}_{R/P}N/PN$. If
$M$ is an $R$-module, $\mathcal{F}(M)$ is the family of its finitely generated submodules and we put 
$\nu(M) = \sup \{ \mu(M')\mid M'\in\mathcal{F}(M)\}$. 

\begin{lemma} \label{L:rank} Let $E$ be a non-zero fp-injective module. If
  $\nu(E)\leq n$ then there exists a non-zero uniserial pure submodule
  $U$ of $E$ such that $\nu(E/U) \leq n-1$. 
\end{lemma}
\emph{Proof.}  By lemma~\ref{L:puruni} there exists a uniserial
pure submodule $U$. Let $M$ be a submodule of $E/U$ generated by
$\{y_1,\ldots,y_p\}$. We assume that $\mu(M)=p$. Let
$x_1,\ldots,x_p\in E$  such that $y_k = x_k+U$ and $F$ be the
submodule of $E$ generated by $x_1,\ldots,x_p$. 

First we assume that
$U$ is finitely generated, so $U=Rx_0$, for some $x_0\in E$. If $F\cap
U = U$ then  $U\subseteq F$ and $U$ is a pure submodule of $F$.  It
follows that the following sequence is exact:
\[0\rightarrow
\frac{U}{PU}\rightarrow\frac{F}{PF}\rightarrow\frac{M}{PM}\rightarrow 0.\] 
So, we have $\mu (M) =\mu (F) - \mu (U)=p-1$. We get a contradiction
since $\mu(M)=p$. Hence $F\cap U \not= U$. Let $N$ be the submodule of
$E$ generated by $x_0,x_1,\ldots,x_p$. Clearly  $Rx_0 = N\cap
U$. It follows that $N/(N\cap U)\simeq M$. Then $\mu (N) =
p+1$. Hence  $p\leq n-1$. 

Now, suppose that $U$ is not finitely generated. So, we may assume
that $(0:x)\ne 0$ for each $x\in E$. Since $U$ is faithful, there exists
$x_0\in U$ such that $(0:x_0) \subset \cap_{i=1}^{i=p}
(0:x_i)$. Clearly $F\cap U\ne U$.  Let
$N$ be the submodule of $E$ generated by $x_0,x_1,\ldots,x_p$. Then $Rx_0$
is a pure submodule of $N$: see proof of \cite[Theorem 2]{War}.  
 Clearly  $Rx_0 = N\cap U$. As above we conclude that $p\leq n-1$. 
 \qed

\bigskip
\emph{Proof of proposition~\ref{P:main}.}~ Let $\mathcal{F}$ be the
family of non-zero ideals $A$ of $R$ for which $R/A$ is not
maximal. If $\mathcal{F}\ne\emptyset$ we set
$J=\cup_{A\in\mathcal{F}}A$. Then $J$ is a prime
  ideal by \cite[Lemma II.6.5]{FuSa}. Since $R$ is archimedean, we
  have either $J=P$ or $\mathcal{F}=\emptyset$. So $R$ is almost 
maximal if and only if
there exists a non-zero proper ideal $I\ne P$ such that $R/I$ is 
maximal. Consequently, we may replace $R$ by $R/rR$ for some nonzero
element $r$ of $P$ and assume that $R$ is an IF-ring by \cite[Theorem
11]{Cou1}. Let $E$ be a non-zero injective indecomposable $R$-module. Since
$\nu(E) \leq n$, by using lemma~\ref{L:rank}, there exists a
pure-composition series of $E$ of length $m\leq n$, 
\[0 = E_0 \subset E_1 \subset \ldots \subset E_m = E\]
such that \   $E_k/E_{k-1}$ \  is uniserial and $fp$-injective, 
$\forall k$, $1\leq k\leq m$.  By lemma~\ref{L:inj}, $U= E/ 
E_{m-1}$ is injective. Since $U$ is faithful by lemma~\ref{L:puruni}
there exists $x\in U$ such
that $\{0\}\subset (0 : x) \subset P$. We set $I = (0: x)$. Let  $V=\{
y\in U \mid I\subseteq (0 : y)\}$. Then $V$ is an injective module over
$R/I$. Assume that there exists $y\in V$ such that $x=sy$ for some
$s\in R$. By \cite[Lemma  2]{Cou1} we get $I\subseteq (0:y)=sI$. As in
\cite[p. 69]{FuSa} let $I^{\sharp}=\{r\in R\mid rI\subset I\}$. By
\cite[Lemma II.4.3]{FuSa} $I^{\sharp}$ is a prime ideal. So, $I^{\sharp}=P$
 because $P$ is the only prime ideal. The equality $sI=I$ implies that
 $s$ is a unit. It follows that $V=Rx\simeq R/I$. So
$R/I$ is self-injective. By \cite[Theorem 2.3]{Kla} $R/I$ is maximal. 
Hence $R$ is almost maximal. \qed

\bigskip
The proof of the following theorem is now complete.

\begin{theorem} Every commutative local ring of bounded module 
type is an almost maximal valuation ring.
\end{theorem}

\bigskip
From proposition~\ref{P:main} we also deduce the following corollary.

\begin{corollary} Let $R$ be an archimedean valuation domain, 
$\widetilde R$ a maximal immediate extension of $R$, $Q$ and $\widetilde 
Q$ their respective fields of fractions. If $[\widetilde Q : Q] <
\infty$, then $R$ is almost maximal. 
\end{corollary}
\emph{Proof.} The conclusion holds since  $\mu (M) \leq [\widetilde Q
: Q]$ for every finitely generated uniform module $M$ by \cite[Theorem
2.2]{Zan1}. \qed

\end{document}